\documentclass[11pt]{article}
\usepackage{graphicx}
\usepackage{amssymb}

\newtheorem{theorem}{Theorem}

\newtheorem{lemma}{Lemma}

\title{Symplectic spreads and permutation polynomials}
\author{
{\em Simeon Ball}
\thanks{This author acknowledges the support of the Ministerio de Ciencia y Tecnologia, Espa{\~n}a.}\\ 
Departament de Matem{\a` a}tica Aplicada IV \\ Universitat Polit{\a` e}cnica de Catalunya \\
 M{\a`o}dul C3, Campus Nord \\  08034 Barcelona, Espanya. \\ {\tt simeon@mat.upc.es} \and
{\em Michael Zieve}\\ Center for Communications Research\\ 805 Bunn Drive \\Princeton, NJ 08540--1966\\
 USA.\\ {\tt zieve@idaccr.org}
\\ }

\date{22 September 2003}
\begin{document}
\maketitle

\begin{abstract}
Every symplectic spread of $PG(3,q)$, or equivalently every ovoid of $Q(4,q)$, is shown to give a
 certain family of permutation polynomials of $GF(q)$ and conversely. This leads to an algebraic
 proof of the existence of the Tits-L{\"u}neburg spread of $W(2^{2h+1})$ and the Ree-Tits spread 
of $W(3^{2h+1})$, as well as to a new family of low-degree permutation polynomials over $GF(3^{2h+1})$.
\end{abstract}

Let $PG(3,q)$ denote the projective space of three dimensions over $GF(q)$. A {\em spread} of $PG(3,q)$
 is a partition of the points of the space into lines. A spread is called {\em symplectic} if every
 line of the spread is totally isotropic with respect to a fixed non-degenerate alternating form.
 Explicitly, the points of $PG(3,q)$ are equivalence classes of nonzero vectors $(x_0,x_1,x_2,x_3)$
 over $GF(q)$ modulo multiplication by $GF(q)^*$. Since all non-degenerate alternating forms on $PG(3,q)$
are equivalent (cf.\ \cite[p.\ 587]{Lang93} or \cite[p.\ 69]{Taylor92}), we may use the form
\begin{equation} \label{form}
((x_0,x_1,x_2,x_3),(y_0,y_1,y_2,y_3))=x_0y_3-x_3y_0-x_1y_2+y_1x_2.
\end{equation}
Then a symplectic spread is a partition of the points of $PG(3,q)$ into lines such that
$(P,Q)=0$ for any points $P,Q$ lying on the same line of the spread.

Symplectic spreads are equivalent to other objects. A symplectic spread is a spread of the
 generalised quadrangle $W(q)$ (sometimes denoted as $Sp(4,q)$), whose points are the points
 of $PG(3,q)$ and whose lines are the totally isotropic lines with respect to a non-degenerate
 alternating form. By the Klein correspondence (see for example \cite{Cameron91},
 \cite[pp. 189]{Taylor92} or \cite{Tits62}), a spread of $W(q)$ gives an ovoid of the
 generalised quadrangle $Q(4,q)$ (sometimes denoted $O(5,q)$) and vice-versa.

Let ${\cal S}$ be a spread of $PG(3,q)$. There are $q^3+q^2+q+1$ points in $PG(3,q)$, and
 each line contains $q+1$ points.
Since $\cal S$ is a partition of the points of $PG(3,q)$ into lines, it contains exactly
 $q^2+1$ lines. The group $PGL(4,q)$ acts transitively on the lines of $PG(3,q)$, so let us
 assume that $\cal S$ contains the line
$l_{\infty}$, which we define as
 $$\langle (0,0,0,1), (0,0,1,0) \rangle.$$  The plane $X_0=0$ contains $l_{\infty}$, so each
 of the other $q^2$ lines of the spread contains precisely one of the $q^2$ points 
$\{ \langle (0,1,x,y)\rangle \ | \ x,y \in GF(q) \}. $ The plane $X_1=0$ also contains $l_{\infty}$,
 so the other $q^2$ lines of the spread are given by two functions $f$ and $g$ such that
 $${\cal S}= l_{\infty} \cup \{ \langle (0,1,x,y), (1,0,f(x,y),g(x,y) \rangle \ | \ x,y \in GF(q) \}.$$
 The spread condition is satisfied if and only if for each $a \in GF(q)$ the plane $X_1=aX_0$ is
 partitioned by the lines of the spread. These planes contain $l_{\infty}$ and meet the other lines
 of ${\cal S}$ in the points
$$\{ \langle (1,a,ax+f(x,y),ay+g(x,y)) \rangle \ | \ x,y \in GF(q) \}.$$ Hence the spread condition
 is satisfied if and only if
$$ (x,y) \mapsto (ax+f(x,y),ay+g(x,y)) $$
is a permutation of $GF(q)^2$ for all $a \in GF(q)$.

We are interested here in symplectic spreads. The line $l_{\infty}$ is totally isotropic with
 respect to the form (\ref{form}).
The other lines of the spread are totally isotropic with respect to the form (\ref{form}) if and
 only if for all $x$ and $y\in GF(q)$ 
$$((0,1,x,y), (1,0,f(x,y),g(x,y))=-y-f(x,y)$$
is zero. Hence $${\cal S}:= l_{\infty} \cup \{ \langle (0,1,x,y), (1,0,-y,g(x,y) \rangle \ |
 \ x,y \in GF(q) \}$$ will be a symplectic spread if and only if
\begin{equation} \label{simplecond}
 (x,y) \mapsto (ax-y,ay+g(x,y)) 
 \end{equation}
is a permutation of $GF(q)^2$ for all $a \in GF(q)$.
Now make the substitution $b=ax-y$ to see that this is equivalent to $x \mapsto a(ax-b)+g(x,ax-b)$
 being a permutation of $GF(q)$ for all $a,b \in GF(q)$, which is equivalent to
$x \mapsto g(x,ax-b) + a^2x$
being a permutation of $GF(q)$ for all $a,b \in GF(q)$.

Although merely an observation, this fact seems not to have been noted before, and as we shall see
 it can be quite useful. So let us formulate this in a theorem.
\begin{theorem} \label{one}
The set of totally isotopic lines 
$$l_{\infty} \cup \{ \langle (0,1,x,y), (1,0,-y,g(x,y)) \rangle \ | \ x,y \in GF(q) \}$$
is a (symplectic) spread if and only if $$x \mapsto g(x,ax-b) + a^2x$$
is a permutation of $GF(q)$ for all $a,b \in GF(q)$.\hfill {\footnotesize $\blacksquare$}
\end{theorem}

Symplectic spreads of $PG(3,q)$ are rare. 
All the known examples are given in Table 1 which comes from \cite{PW00}.
In particular, the \emph{regular} spreads are those for which the polynomial
$g(x,y)$ has total degree 1.
The main result in \cite{BGS03} implies that when $q$ is prime, symplectic spreads of $PG(3,q)$ are regular.

\begin{table} \label{examples}
\begin{center}
\begin{tabular}{|l|c|c|c|}
\hline
  name & $g(x,y)$ & $q$ & restrictions \\
\hline
  regular & $-nx$ & odd & $n$ non-square\\
  Kantor \cite{Kantor82} & $-nx^{\alpha}$ & odd & $n$ non-square, $\alpha|q$ \\
  Thas-Payne \cite{TP94} & $-nx-(n^{-1}x)^{1/9}-y^{1/3}$ & $3^h$ & $n$ non-square, $h>2$\\
  Penttila-Williams \cite{PW00} & $-x^9-y^{81}$ & $3^5$ &  \\
  Ree-Tits slice \cite{Kantor82} & $-x^{2\alpha+3}-y^{\alpha}$ & $3^{2h+1}$
   & $\alpha=\sqrt{3q}$ \\  \hline

  regular & $cx+y$ & even & $Tr_{q \rightarrow 2}(c)=1$ \\
  Tits-L{\"u}neburg  \cite{Tits62} & $x^{\alpha+1}+y^{\alpha}$ & $2^{2h+1}$ &
$\alpha=\sqrt{2q}$ \\ \hline
\end{tabular}
\end{center}
\caption{The known examples of symplectic spreads of $PG(3,q)$}
\end{table}

Note that from any of the examples in the table we could make many other equivalent symplectic
 spreads and that the function $g(x,y)$ will not in general have such a nice form. For instance,
 all the examples in the table give spreads
${\cal S}$ that contain the line $l$
$$ \langle (0,1,0,0),(1,0,0,0) \rangle.$$
The linear map $\tau$ that switches $X_0$ and $X_3$ and switches $X_1$ and $X_2$ preserves the
 form (\ref{form}) but swiches $l_{\infty}$ and $l$. The other $q^2-1$ lines in ${\cal S}$ are
 mapped to the lines
$$\{ \langle (y,x,1,0), (g(x,y),-y,0,1) \rangle \ | \ x,y \in GF(q), \ (x,y) \neq (0,0) \}$$
by $\tau$.
Writing these lines as the spans of their points on the planes $X_0=0$ and $X_1=0$, these lines are
$$\{ \langle (0,1,u,v), (1,0,-v,{-v x \over y}) \rangle \ | \ x,y \in GF(q), \ (x,y) \neq (0,0)\},$$
where $$u={g(x,y) \over xg(x,y)+y^2}$$ and $$v={-y \over xg(x,y)+y^2}.$$ 
(When $y=0$ we interpret $-vx/y$ to be $1/g(x,0)$.)
Now one would have to calculate $-vx/y$ in terms of $u$ and $v$ to deduce the function $g(x,y)$ for
 the spread $\tau({\cal S})$.
For an explicit example of this, consider the Kantor spread ${\cal S}$
over
$GF(27)$ with $g(x,y)=-nx^3$, where $n^3-n=-1$.
The function $g(u,v)$ for $\tau({\cal S})$ is
$$
nu^{21}v^4 + n^8u^{19}v^{18} + n^2u^{17}v^6 + n^4u^9v^{10} + n^{18}u^5v^{12} + 
    n^{12}u^3.
$$

A polynomial $h$ in one variable over $GF(q)$ is called {\em additive}
if $h(x+u)=h(x)+h(u)$ for all $x,u\in GF(q)$.  In case $g(x,y)=h_1(x)+h_2(y)$
with $h_1$ and $h_2$ being additive polynomials,
the symplectic spread corresponds to a translation ovoid of $Q(4,q)$, which in
turn comes from a semifield flock of the quadratic cone in $PG(3,q)$.
This has been the subject matter of a number of articles, see for example
\cite{BL94}, \cite{BBL03} or \cite{Lavrauw01}. The classification of such
examples is an open problem whose solution would be of much interest.
The partial classification in \cite{BBL03} implies that if there are any
further examples over $GF(p^h)$ then $p<4h^2-8h+2$. 
Theorem~\ref{one} in this case reads: The polynomial $g(x,y)=h_1(x)+h_2(y)$ 
will give a symplectic spread if and only if
$h_1(x)+h_2(ax)+a^2x$
is a permutation polynomial for all $a\in GF(q)$, or equivalently
$h_1(x)+h_2(ax)+a^2x$ has no zeros in $GF(q)^*$ for all $a\in GF(q)$.

The two examples where $g(x,y)$ is not of this form are the Tits-L{\"u}neburg spread and
 the Ree-Tits spread.

Let us first check the Tits-L{\"u}neburg example, where $\alpha=\sqrt{2q}$. In this case
$$g(x,ax-b) + a^2x=x^{\alpha+1}+(ax)^{\alpha}-b^{\alpha}+a^2x.$$
So we should have that $x^{\alpha+1}+(ax)^{\alpha}+a^2x$ is a permutation polynomial
for all $a\in GF(q)$, which is easy to see since this polynomial is
 $(x+a^{\alpha})^{\alpha+1}+a^{\alpha+2}$. Note that composing permutation polynomials
 with permutation polynomials gives permutation polynomials so it is enough to check
 that $x^{\alpha+1}$ is a permutation polynomial, which it is since 
$(2^{h+1}+1,2^{2h+1}-1)=1$.

%

Now we come to the interesting Ree-Tits slice example, $g(x,y)=-x^{2\alpha+3}-y^{\alpha}$
 where $q=3^{2h+1}$ and $\alpha=\sqrt{3q}$.
This spread was discovered by Kantor \cite{Kantor82} as an ovoid of $Q(4,q)$. It is the
 slice of the Ree-Tits ovoid of $Q(6,q)$. It provides us with an interesting class of
 permutation polynomials, namely, the polynomials $f_a(x):=b^{\alpha}-(g(x,a x-b)+a^2x)$, 
$$f_a(x)=x^{2\alpha+3}+(ax)^{\alpha} - a^2 x.$$
The polynomial $f_a$ is remarkable in that it is a permutation
polynomial over $GF(q)$ whose degree is approximately $\sqrt{q}$.
There are only a handful of known permutation polynomials with
such a low degree.  The bulk of these examples are \emph{exceptional polynomials},
namely polynomials over $GF(q)$ which permute $GF(q^n)$ for infinitely many values $n$.
However, we will show below that $f_a$ is not exceptional, so long as $\alpha>3$ and $a\ne 0$.
There are also some non-exceptional permutation polynomials of degree approximately $\sqrt{q}$
in case $q$ is a square or a power of 2.  However, our example is the first for which $q$ is
an odd nonsquare.

It follows from \cite{Kantor82} and Theorem~\ref{one} that $f_a$ is a permutation polynomial.
 Conversely we now give a direct proof that $f_a$ is a permutation polynomial, which
 (along with Theorem~\ref{one}) gives a new proof that the Ree-Tits examples are in fact
 symplectic spreads.
Our proof that $f_a$ is a permutation polynomial uses the method of Hans Dobbertin
\cite{Dobbertin03}.

\begin{theorem} \label{pps}
Let $q=3^{2h+1}$ and let $\alpha=\sqrt{3q}$. For all $a\in  GF(q)$
the polynomial $f_a(x):=x^{2\alpha+3}+(a x)^{\alpha} - a^2 x$ is a
permutation polynomial over $GF(q)$.
\end{theorem}

{\it Proof.}
If $f_a(x)$ is a permutation polynomial then so is
$\zeta^{2\alpha+3}f_a(x/\zeta)$ for any $\zeta\in GF(q)^*$, and
the latter polynomial equals $f_{a\zeta^{\alpha+1}}(x)$.
Since $(\alpha+1,q-1)=2$, it follows that if $f_a$ is
a permutation polynomial then so is $f_{a\zeta^2}$ for any
$\zeta\in GF(q)^*$.  Thus it suffices to verify the theorem
for a single nonzero square $a$, a single nonsquare $a$,
and the value $a=0$ (in which case the theorem is trivial).
Since $-1$ is a non-square in $GF(3^{2h+1})$ we can assume from now on that $a^2=1$.

Suppose that $f_a$ is not a permutation polynomial.  Let $x,y$ be distinct
elements of $GF(q)$ such that $f_a(x)=f_a(y)=d$. The equations $f_a(x)=d$ and
$f_a(x)^{\alpha}=d^{\alpha}$ give
\begin{eqnarray}
\label{1} x^{2{\alpha}+3}+a x^{\alpha} -  x = d \\
\label{2} x^{6+3\alpha}+a x^3 -  x^{\alpha} = d^{\alpha}.
\end{eqnarray}
By viewing these equations as low-degree polynomials in $x^{\alpha}$ whose
coefficients are low-degree polynomials in $x$, we can solve for $x^{\alpha}$
as a low-degree rational function in $x$.  Namely, multiplying (\ref{1})
by $x^{\alpha+3}$ and then subtracting (\ref{2}) gives
\begin{equation}
\label{3}a x^{2\alpha+3} - x^{\alpha+4} - a x^3 +
 x^{\alpha} = d x^{\alpha+3} - d^{\alpha};
\end{equation}
multiplying (\ref{1}) by $a$ and subtracting (\ref{3}) gives
\begin{equation}
\label{4}x^{\alpha} (x^4 + d x^3) = a x + d a
 - a x^3 + d^{\alpha}.
\end{equation}
This expresses $x^{\alpha}$ as a low-degree rational function in $x$, so
long as $x\notin\{0,-d \}$.  For later use we record this
equation in the form $F(x^{\alpha},x)=0$ where
$$
F(T,U):= U^4 T + d U^3 T - a U - d a
 + a U^3 - d^{\alpha}.
$$
Note that $x$ and $y$ are not both in $\{0,-d \}$, for if so
then $d=f_a(0)=0$ so $x=y=0$, contradiction.  Thus, by swapping
$x$ and $y$ if necessary, we may assume $x\notin\{0,-d \}$.

Solving for $x^{\alpha}$ in (\ref{4}) and substituting into (\ref{1}) gives a
low-degree polynomial satisfied by $x$:
$$
(a x + da - a x^3 + d^{\alpha})^2  + 
 a ( x + d) (a x + d a -
a x^3 + d^{\alpha}) = ( x^2 + d x)^3.
$$
By expanding this equation we get
\begin{equation} \label{5} 
(d^{\alpha}a-d^3)x^3- x^2+d  x-d^2 +d^{2\alpha}=0.
\end{equation}
Next we handle the cases $y=0$ and $y=-d$.
If $y=0$ then $d=f_a(y)=0$ and (\ref{5}) implies $x=0$, contradiction. 
If $y=-d$ then the analogue of
(\ref{4}) with $y$ in place of $x$ says that $d^{\alpha}=-a d^3$, so
$d^{2\alpha-6}=1$ and since $(q-1,2\alpha-6)=2$ that $d^2=1$ and hence $a=-1$.
Then equation (\ref{5}) simplifies to $dx(x+d)^2=0$. Since $d=0$ implies $x=0$ we have
 $x\in\{0,-d \}$, again a contradiction.

Hence we may assume $y\notin\{0,-d\}$, and moreover we may assume
$d\ne 0$ and $d^3\ne  -d^{\alpha}a$. We can also assume that $d^3\ne  d^{\alpha}a$.
For, if $d^3= d^{\alpha}a$ then $d^{2\alpha-6}=1$ and again 
since $(q-1,2\alpha-6)=2$ that $d^2=1$ and hence $a=1$. Then equation (\ref{5})
 simplifies to $x(x+d)=0$, a contradiction.

In particular, equation (\ref{5}) remains valid if we
substitute $y$ for $x$.
Thus $x$ and $y$ are roots of the polynomial
\begin{equation}
\psi(t):=(d^{\alpha}a-d^3)t^3-t^2+d t-d^2+d^{2\alpha}.
\end{equation}
We express the roots of $\psi(t)$ in terms of $x$.
Since $\psi(x)=0$, we know that $t-x$ is a factor of $\psi(t)$: in fact,
writing $A:=d^{\alpha}a-d^3$, we have
\begin{equation}
\psi(t)/(t-x) = At^2 + (Ax-1)t + (Ax^2+d-x).
\end{equation}
The discriminant of this quadratic polynomial is
$$\delta:=(Ax-1)^2-A(Ax^2+d-x)=1-A(x+d).$$
If $\delta=0$ then $x=-d+1/A$ and $y=(Ax-1)/A=-d$ which we have already excluded,
 so assume from now on that $\delta\neq 0$.

Substituting $d=x^{2\alpha+3}+a x^{\alpha}-x$ we find that 
$$
A=-x^{6\alpha+9}+ax^{3\alpha+6}-ax^{3\alpha}-ax^{\alpha}-x^3
$$
and
$$\delta=(x^{4\alpha+6}-a x^{3\alpha+3}+x^{2\alpha}+ax^{\alpha+3}-1)^2.$$
Thus putting $\sqrt{\delta}= x^{4\alpha+6}-a x^{3\alpha+3}+x^{2\alpha}+ax^{\alpha+3}-1$,
 we can write the roots of $\psi(t)/(t-x)$ as
$$
y_ 1:= x-(\sqrt{\delta}+1)/A={x^{3\alpha+4}+ax^{2\alpha+1}+x^{\alpha}+ax \over
 x^{3\alpha+3}+ax^{2\alpha}+a}
$$
and
$$
y_2:= x+(\sqrt{\delta}-1)/A={x^{3\alpha+7}-ax^{2\alpha+4}-x^{\alpha+3}+x^{\alpha+1}+ax^4-a
 \over x^{3\alpha+6}-ax^{2\alpha+3}+x^{\alpha}+ax^3}.
$$
Now one can verify that $F(y_2^{\alpha},y_1)=0$ and $F(y_1^{\alpha},y_2)=0$.
But we know that $F(y^{\alpha},y)=0$ and $y\in\{y_1,y_2 \}$.  Since $y_1\ne y_2$,
this implies $F(T,y)=0$ has more than one root.  But this is a linear
polynomial in $T$, a contradiction. \hfill {\footnotesize $\blacksquare$}

Recall that a polynomial $f$ over $GF(q)$ is called \emph{exceptional} if it permutes
$GF(q^n)$ for infinitely many $n$.  We now show that, except in some special cases,
$f_a$ is not exceptional.  Our proof relies on the classification of monodromy groups
of indecomposable exceptional polynomials, due to Fried, Guralnick, and Saxl \cite{FGS}.
A polynomial is {\em indecomposable} if it is not the
composition of two polynomials of lower degree.

\begin{lemma}
When $\alpha>3$ and $a\ne 0$,  $f_a(x)$ is indecomposable.
\end{lemma}

{\it Proof.} 
The derivative of $f_a$ is $f_a'=-a^2$, which is a
nonzero constant. If $f_a(x)=g(h(x))$ then $-a^2=f_a'(x)=g'(h(x))h'(x)$,
so both $g'$ and $h'$ are nonzero constants. Thus
$g(x)=u(x^3)+cx$  and  $h(x)=v(x^3)+dx$  for some polynomials $u$ and $v$ and nonzero
 constants $c$ and $d$. Since the degree of $f_a$ is not divisible by $9$,
either $g$ or $h$ has degree not divisible by $3$, and hence must have degree 1.
Thus $f_a$ is indecomposable.\hfill {\footnotesize $\blacksquare$}

\begin{theorem}
When $\alpha > 3$ and $a\ne 0$,  $f_a(x)$ is not exceptional.
\end{theorem}

{\it Proof.} 
This follows directly from the preceeding lemma and \cite[Theorems 13.6 and 14.1]{FGS},
 according to which there is no
indecomposable exceptional polynomial of degree $2\alpha+3$ over a finite
field of characteristic $3$.\hfill {\footnotesize $\blacksquare$}

In all the examples in Table 1 the polynomial $g$ is of the form $g(x,y)=h_1(x)+h_2(y)$.
 In Glynn \cite{Glynn84} such a polynomial $g(x,y)$ with this property is called {\em separable}.
 Every known example of a symplectic spread of $PG(3,q)$ is equivalent to a symplectic
 spread with $g(x,y)$ separable. In the examples not only is the polynomial $g(x,y)=h_1(x)+h_2(y)$
 separable but $h_2(y)=Cy^{\sigma}$, where $y \mapsto y^{\sigma}$ is an automorphism of $GF(q)$. 
We can classify these examples in the case when $q$ is even using Glynn's Hering classification of
 inversive planes \cite{Glynn84}.

\begin{theorem}
Let $q$ be even. If $g(x,y)=h_1(x)+Cy^{\sigma}$ is a separable polynomial that gives a symplectic 
spread of $PG(3,q)$ then the spread is either a regular spread or a Tits-L{\"u}neburg spread.
\end{theorem}

{\it Proof.}
If $C=0$ then Theorem~\ref{one} implies
$h_1(x)+a^2x$ is a permutation polynomial for all $a \in GF(q)$.
Let $x$ and $y$ be distinct elements of $GF(q)$, and put $d=(h_1(x)+h_1(y))/(x+y)$.
Then $h_1(x)+dx=h_1(y)+dy$, so the polynomial $h_1(x)+dx$ is not a permutation polynomial,
 a contradiction.

Now assume that $C \neq 0$. Put $z=h_1(x)+Cy^{\sigma}$ and rewrite this as
 $y=C^{-1}z^{1/\sigma}-C^{-1}h_1(x)^{1/\sigma}$. Define the function
 $s(x,z):=C^{-1}z^{1/\sigma}-C^{-1}h_1(x)^{1/\sigma}$. Then $g(x,y)=z$ if and only if $s(x,z)=y$.
 We have already seen in equation (\ref{simplecond}) that $g(x,y)$ will give a symplectic
 spread if and only if 
$$ (x,y) \mapsto (ax-y,ay+g(x,y)) $$
is a permuatation of $GF(q)^2$. This is equivalent to the condition that for all $(x,y) \neq (u,v)$
 $$(ax-y,ay+g(x,y)) \neq (au-v,av+g(u,v))$$ for all $a \in GF(q)$. If these pairs were equal then
 eliminating $a$ this gives the condition that for all $(x,y) \neq (u,v)$
$$(y-v)^2+(x-u)(g(x,y)-g(u,v))\neq 0.$$
Now put $z=g(x,y)$ so that $s(x,z)=y$, and put $w=g(u,v)$ so that $s(u,w)=v$.
Then we have that  $(x,z) \neq (u,w)$
$$(s(x,z)-s(u,w))^2+(x-u)(z-w) \neq 0.$$
When $q$ is even this is exactly the polynomial condition on such a polynomial $s(x,z)$ that
 Glynn studies in \cite{Glynn84} and that he classifies as coming from either a regular spread
 or a Tits-L{\"u}neburg spread.
\hfill {\footnotesize $\blacksquare$}

When $q$ is odd we can use Thas' classification of flocks of the quadratic cone in $PG(3,q)$
 whose planes are incident with a common point from \cite{Thas87} to prove the following theorem. 
We realise that for many readers familiar with flocks and semifield flocks the next theorem is
 immediate, but we include a proof for those readers who may not be.

\begin{theorem}
Let $q$ be odd. If $g(x,y)=h_1(x)$ is a separable polynomial that gives a symplectic spread
 of $PG(3,q)$ then the spread is either a regular spread or a Kantor spread.
\end{theorem}

{\it Proof.}

Consider the set of $q$ planes of $PG(3,q)$
$$ \{ X_0 + h_1(x)X_1+xX_3=0 \ |\ x \in GF(q) \}.$$  
We claim that any two of these planes intersect in a line which is disjoint
from the degenerate quadric $X_1X_3=X_2^2$. Indeed take two planes coordinatised by $x$ and $y$,
 $x\neq y$. Then the points in their intersection $(z_0,z_1,z_2,z_3)$ satisfy
 $(h_1(x)-h_1(y))z_1+(x-y)z_3=0$, and the points which also 
lie on the degenerate quadric satisfy $$(h_1(x)-h_1(y))z_1^2+(x-y)z_2^2=0.$$ If $z_1 \neq 0$ then
 $h_1(x)+(z_2/z_1)^2x$ is not a permutation polynomial, a contradiction. If $z_1=0$ then $z_2=0$
 and $z_0=-xz_3=-yz_3$. But $x\neq y$ implies that $z_3=0$ and $z_0=0$ which is nonsense.
We have shown that the set of planes form a flock of the quadratic cone in $PG(3,q)$. Moreover
 all these planes are incident with $(0,0,1,0)$. By a theorem of Thas \cite{Thas87} this flock
 is either linear or of Kantor type.
In other words, the spread is either regular or Kantor.
\hfill {\footnotesize $\blacksquare$}

In general the permutation polynomial condition from Theorem~\ref{one} requires the existence
 of $q^2$ permutation polynomials, one for each pair $(a,b) \in GF(q)^2$. If
 $g(x,y)=h_1(x)+h_2(y)$ and $h_2(y)$ is additive then Theorem~\ref{one} simplifies to: The
 polynomial $g(x,y)=h_1(x)+h_2(y)$ 
will give a symplectic spread if and only if
$f_a(x):=h_1(x)+h_2(ax)+a^2x$
is a permutation polynomial
for all $a\in GF(q)$.
This condition only requires the existence of $q$ permutation polynomials. Moreover as we saw
 in the proof of Theorem~\ref{pps}, if the non-zero terms in $h_1$ and $h_2$ have suitable
 degrees, many of these permutation polynomials may be equivalent.

\begin{table} \label{classes}
\begin{center}
\begin{tabular}{|l|c|c|c|}
\hline
  name & $q$ & $\Delta$ & $(q-1,\Delta d)+1$ \\
\hline
  regular & odd & $1$ & $1$\\
  Kantor & odd & $1$ & $(q-1,(\alpha-1)/2)+1$ \\
  Thas-Payne & $3^h$ & . & $q$ \\
  Penttila-Williams & $3^5$ &  $11$ & $23$   \\
  Ree-Tits slice  &   $3^{2h+1}$ & $1$ & $3$\\
  \hline 
  regular &  even & $1$ & $1$ \\
  Tits-L{\"u}neburg  & $2^{2h+1}$ &
$1$ & $2$ \\ \hline
\end{tabular}
\end{center}
\caption{The class $\Delta$ of the known examples of symplectic spreads of $PG(3,q)$}
\end{table}

Let us investigate this further. We define a set of polynomials
$\{f_a(x) \ |æ\ a \in GF(q) \}$ to be of {\it class $\Delta$} if there exists a $t$ and $d$
 such that $$f_a(bx)=b^tf_{ab^d}(x)$$ for all $b^{q-1/\Delta}=1$ and $a$ and $x \in GF(q)$.
 Now we can lessen the condition in Theorem~\ref{one} for $\Delta < q-1$.
\begin{theorem}
Let the set of $q$ polynomials $\{f_a(x) \ |æ\ a \in GF(q) \}$, where $f_a(x)=h_1(x)+h_2(ax)+a^2x$
 and $h_2$ is additive, be of class $\Delta$. The $f_a$
is a permutation polynomial for all $a\in GF(q)$ if and only if $f_a$
is a permutation polynomial for $a=0$ and $a=\varepsilon ^r$, for all $1 \leq r < (q-1,\Delta d)$,
 where $\varepsilon$ is a fixed primitive element.
\end{theorem}

{\it Proof.}
Write $a=\varepsilon^{n_1(q-1,\Delta d)+n_0}$ where $n_0< (q-1,\Delta d)$. Now 
choose $b$ such that $b^d= \varepsilon^{-n_1(q-1,\Delta d)}$. \hfill {\footnotesize $\blacksquare$}

In Table 2 we have listed the class for the known examples and the quantity $(q-1,\Delta d)+1$,
 the number of permutation polynomials that need to be checked in each case. Inspired by this
 table we used the mathematical package GAP to look at polynomials over $GF(q)$, $q=p^h$, of the
 form $g(x,y)=Dx^t+Cy^\sigma$ for all $\sigma$ a power of $p$ and $D$ and $C$ elements of $GF(q)$
 where the corresponding set of polynomials 
$\{f_a(x) \ |æ\ a \in GF(q) \}$ is of class $\Delta$ with $\Delta$ small. An exhaustive search was
 carried out for $\Delta \leq 23$ and $q\leq 67^2=4489$, $\Delta=2$ and $q <3^8=6561$, $\Delta=1$
 and $q<3^9=19683$.
No new examples of symplectic spreads were found.

\bibliographystyle{plain}

\end{document}